\documentclass[11pt,a4paper]{article}
\usepackage{amssymb,amsmath}

\newtheorem{theorem}{Theorem}[section]
\newtheorem{corollary}[theorem]{Corollary}

\newtheorem{lemma}[theorem]{Lemma}
\newtheorem{proposition}[theorem]{Proposition}
\newtheorem{remark}[theorem]{Remark}

\newtheorem{example}[theorem]{Example}

\newenvironment{proof}{\begin{trivlist}\item[]{\it
Proof.}}{\hfill$\square$\end{trivlist}}

\def\field{k}
\def\cov{{\mathrm{Cov}}}
\def\rcov{{\mathrm{RCov}}}
\def\iid{{\mathrm{id}}}
\def\mor{{\mathrm{Mor}}}

\def\gl{{\mathrm{GL}}}
\def\lll{l}

\def\charac{{\mathrm{char}}}
\def\rto{\to}
\def\dif{{\mathrm{d}}}
\def\dom{{\mathrm{dom}}}
\def\spec{{\mathrm{Spec}}} 
\def\pr{{\mathrm{pr}}}
\def\rat{{\mathrm{Rat}}}
\def\imm{{\mathrm{Im}}}

\begin{document}
\title{Covariants and the no-name lemma} 
\author{M\'aty\'as Domokos
\thanks{Partially supported by 
OTKA NK72523 and K61116. } 
\\ 
\\
{\small R\'enyi Institute of Mathematics, Hungarian Academy of 
Sciences,} 
\\ {\small P.O. Box 127, 1364 Budapest, Hungary,} {\small E-mail: domokos@renyi.hu } 
}
\date{}
\maketitle 
\begin{abstract}  
A close connection between the no-name lemma 
(concerning algebraic groups acting on vector bundles) 
and the existence of sufficiently many independent rational covariants is 
pointed out. In particular, this leads to a new natural 
proof of the no-name lemma. 
For linearly reductive groups, the approach has a refined variant based on integral covariants. This fits into the usual context of invariant theory, and yields a version of the no-name lemma that has a constructive nature. 

\end{abstract}



\section{Introduction}\label{sec:intro} 

The so-called rationality problem asks when $\field(X)^G$, the field of invariant rational functions is purely transcendental over $\field$ (an algebraically closed base field), where $G$ is a group acting morphically on the irreducible algebraic variety $X$. 
One of the basic principles in studying this question is the following statement (for simplicity, we state it for the case $\charac(\field)=0$): 
if $G$ acts generically freely on $X$, and $W$ is a finite dimensional $G$-module, then $\field(X\times W)^G$ is purely transcendental over $\field(X)^G$.   This statement (or some variant of it) has gone into the literature as the "no-name lemma" (see Remark~\ref{remark:precise} for references). 

In this note we point out that the no-name lemma is closely related to the question about the number of generically independent (rational) covariants from $X$ to $W$. 
In particular, we show that the no-name lemma follows from the following fact: 
if the stabilizer $G_x$ of a general $x\in X$ acts trivially on $W$, then there are $\dim(W)$ generically independent rational covariants from $X$ to $W$. 
This latter fact was proved by Reichstein \cite{reichstein}
for generically free actions in characteristic zero; we extend it to positive characteristic and the above weaker and necessary condition on stabilizers.   The picture is summarized in Theorem~\ref{thm:main}, stating the equivalence of various properties expressed in terms of stabilizers, rational covariants, birational isomorphisms, and generators of invariant fields, respectively.  

A benefit from paying attention to covariants is that it leads directly to a stronger form of the no-name lemma, stated before only for finite groups. 
Namely, if $G_x$ acts trivially on $W$ for a general $x\in X$ (and the action 
is generically separable), then there is a birational isomorphism between 
$X\times W$ and $X\times\field^{\dim(W)}$ ($G$ acting trivially on the vector space $\field^{\dim(W)}$) that restricts to a $\field$-linear isomorphism $W\to\field^{\dim(W)}$ for a general $x\in X$. 
(For finite groups, this is essentially Speiser's lemma applied for function fields of $G$-varieties.)   

The main advantage of this approach is that it has a finer variant for linearly reductive groups, that fits nicely into the traditions of invariant theory. We find exact conditions ensuring the supply of integral covariants to produce $G$-equivariant isomorphisms over invariant affine open sets defined by the non-vanishing of some relative invariant, see Corollary~\ref{cor:2reductive} for a precise statement. 
Moreover, since there are good algorithms to compute generators of modules of covariants for linearly reductive groups, this can be viewed as a constructive version of the no-name lemma (see Remark~\ref{remark:constructive}). 
There are also a couple of interesting examples where the covariants needed 
come naturally with the situation, see Section~\ref{sec:examples}. 


\section{Rational covariants} \label{sec:mainlemma} 

Let $G$ be a linear algebraic group over an algebraically closed field 
$\field$ of arbitrary characteristic, denote by $G_0$ the connected component of the identity. 
By a {\it $G$-variety} $X$ we mean an algebraic variety with a $G$-action such that the action map $G\times X\to X$ is a morphism of algebraic varieties. 
Write $\field[X]$ for the ring of regular functions on $X$, and when $X$ is irreducible, write $\field(X)$ for the field of rational functions on $X$. 
As usual, $\field[X]^G$ is the subring of $G$-invariants, and 
$\field(X)^G$ is the subfield of $G$-invariants.  
When $X$ is affine, $\field[X]$ is called the coordinate ring of $X$, and 
$\field(X)$ is the field of fractions of $\field[X]$ in this case. 
By a {\it $G$-module} we mean a $G$-variety $V$ which is a finite dimensional vector space, on which $G$ acts via linear transformations. 
We shall write $\field^d$ for the $d$-dimensional vector space endowed with the trivial $G$-action. 
We say that some property {\it holds for a general $x\in X$} if it holds for all $x\in U$, where $U$ is a Zariski dense open subset of $X$. 

Let $X$ be an irreducible $G$-variety and $W$ a $G$-module. 
A {\it covariant} $F:X\to W$ is a $G$-equivariant morphism of algebraic varieties. Write $\cov_G(X,W)$ for the set of covariants from $X$ to $W$; 
this is contained in $\rcov_G(X,W)$, the set of {\it rational covariants} 
(i.e. rational $G$-equivariant maps). Note that $\rcov_G(X,W)$ is naturally a vector space over $\field(X)^G$, whereas $\cov_G(V,W)$ is a 
$\field[X]^G$-submodule. 

We say that the (rational) covariants $F_1,\ldots,F_d$ are {\it generically independent}, 
if $F_1(x),\ldots,F_d(x)$ are linearly independent vectors in $W$ for some 
$x\in X$. Note that in this case $F_1(x),\ldots,F_d(x)$ are linearly independent for a general $x\in X$. 

For $f\in\field[X]$ denote by $X_f$ the Zariski open subset 
$\{x\in X\mid f(x)\neq 0\}$. 
Recall that $f\in\field[X]$ is a {\it relative invariant} if 
$g\cdot f=\theta(g)f$ for some character $\theta:G\to\field^{\times}$. 
In this case $X_f$ is a $G$-stable subset of $X$. 
Let us introduce the following ad hoc terminology: we say that the $G$-varieties 
$X\times Y$ and $X\times Z$ {\it are (birationally) isomorphic over} $X$ if there is a 
$G$-equivariant (birational) isomorphism between them that commutes with the projection onto $X$.

We start with an elementary lemma: 

\begin{lemma}\label{lemma:indep-cov} 
Let $X$ be an irreducible $G$-variety, $W$ a $d$-dimensional $G$-module, 
and $F_1,\ldots,F_d\in\cov_G(X,W)$ generically independent covariants. Then there exists a non-zero relative invariant $f\in\field[X]$ such that the $G$-varieties  $X_f\times W$ and  $X_f\times\field^d$ ($G$ acting trivially on the vector space $\field^d$) 
are isomorphic over $X_f$ (via an isomorphism constructed explicitly in the proof). 
\end{lemma} 

\begin{proof} Suppose $F_1,\ldots,F_d\in\cov_G(X,W)$ are generically independent.  
Let $e_1,\ldots,e_d$ be a basis in $W$, and $\varepsilon_1,\ldots,\varepsilon_d$ the corresponding dual basis in $W^*$. 
Then $F_j(x)=\sum_{i=1}^dF_{ij}(x)e_i$ for some $F_{ij}\in\field[X]$. 
Write $F$ for the $d\times d$ matrix whose $(i,j)$-entry is $F_{ij}$. 
The covariance of the $F_j$ and multiplicativity of the determinant imply that 
$f:=\det(F_{ij})_{d\times d}$ is a relative invariant in $\field[X]$ 
of weight $g\mapsto \det(g_W)^{-1}$, where $g_W$ is the matrix of $g$ acting on $W$ with respect to the chosen basis. Moreover, $X_f$ is the locus of $x\in X$ where 
the $F_1(x),\ldots,F_d(x)$ span $W$.  
Define 
$\Phi=(\Phi_1,\ldots,\Phi_d):X_f\times W\to\field^d$ by 
\begin{equation}\label{eq:1phidef}
w=\sum_{i=1}^d\Phi_i(x,w)F_i(x).
\end{equation}
Applying $g\in G$ to both sides, and taking into account the linearity 
of the action of $G$ on $W$ one gets that the maps 
$\Phi_i:X_f\times W\to\field$ are $G$-invariant, so the morphism 
$(\iid_{X_f}\times\Phi): X_f\times W\to X_f\times \field^d$ 
is indeed $G$-equivariant. 

In terms of coordinates, $\Phi_i=\sum_{j=1}^d(F^{-1})_{ij}\varepsilon_j$, where $(F^{-1})_{ij}$ is the $(i,j)$-entry 
of the inverse of $F$. 
This shows that $(\iid_{X_f}\times \Phi):X_f\times W\to X_f\times\field^d$ 
is a morphism of algebraic varieties, and $\Phi(x,-):W\to\field^d$ is a 
$\field$-linear isomorphism for all $x\in X_f$. Moreover, formula 
(\ref{eq:1phidef}) shows that $\iid_{X_f}\times\Phi$ is in fact an isomorphism, with inverse sending $(x,a)\in X_f\times\field^d$ to 
$(x,\sum_{i=1}^da_iF_i(x))\in X_f\times W$. 
\end{proof} 

Lemma~\ref{lemma:indep-cov} has the following converse: 

\begin{lemma}\label{lemma:converse-indep-cov} 
Let $Y$ be an irreducible $G$-variety, $W$ a $d$-dimensional $G$-module, and suppose that the $G$-varieties 
$Y\times W$ and $Y\times \field^d$ are isomorphic over $Y$. 
Then there is a $G$-equivariant isomorphism over $Y$ between them which restricts to a $\field$-linear isomorphism 
$\{x\}\times W\cong W\to \field^d\cong\{x\}\times \field^d$ for all $x\in Y$, and 
there exist $d$ covariants $F_i:Y\to W$ $(i=1,\ldots,d)$ such that $F_1(x),\ldots,F_d(x)$ are linearly independent for all 
$x\in Y$. 
\end{lemma} 

\begin{proof} Let $(x,w)\mapsto (x,\Phi(x,w))$ be a $G$-isomorphism $Y\times W\to Y\times \field^d$. 
Consider the coordinate functions $\Phi_i$ $(i=1,\ldots,d)$ of $\Phi$. 
Then $\Phi_i\in\field[Y\times W]^G$. View $\field[Y\times W]$ as a polynomial ring in the variables $\varepsilon_1,\ldots,\varepsilon_d$ 
(a basis of $W^*$) with coefficients in $\field[Y]$. 
The linear component of $\Phi_i$ is $\sum_{j=1}^d\Phi_{ij}\varepsilon_j$ with 
$\Phi_{ij}\in\field[Y]$. Since the action of $G$ is homogeneous, 
we have that $\sum_{j=1}^d\Phi_{ij}\varepsilon_j$ is $G$-invariant. 
Moreover, for all $x\in Y$, the matrix $(\Phi_{ij}(x))_{d\times d}$ is invertible, being the matrix of the differential at zero of the isomorphism $\Phi(x,-):W\to\field^d$, so $\det(\Phi_{ij})_{d\times d}$ is a unit in $\field[Y]$. 
The $d$ desired covariants are $x\mapsto\sum_{i=1}^dF_{ij}(x)e_i$, $j=1,\ldots,d$, 
where $F_{ij}(x)$ is the $(i,j)$-entry of the inverse of the matrix 
$(\Phi_{st}(x))_{d\times d}$, and $e_1,\ldots,e_d$ is the basis dual to 
$\varepsilon_1,\ldots,\varepsilon_d$. 
Indeed, the $G$-invariance of the $\Phi_i$ can be expressed as the matrix equality 
\begin{equation}\label{eq:2matrixeq} 
(g\Phi_{ij})_{d\times d}=(\Phi_{ij})_{d\times d}\cdot g_W,  
\end{equation} 
where $g_W$ is the matrix of $g\in G$ acting on $W$ with respect to a basis 
$e_1,\ldots,e_d$ dual to $\varepsilon_i$.
Formula (\ref{eq:2matrixeq}) shows that 
$F_j:x\mapsto\sum_{i=1}^dF_{ij}(x)e_j$ 
is a covariant for $j=1,\ldots,d$. 
Moreover, since $\det(F_{ij})_{d\times d}$ is a unit in $\field[Y]$, the $F_j(x)$ $(j=1,\ldots,d)$ are 
linearly independent for all $x\in Y$. 
\end{proof} 

It is a bit less obvious, but turns out from Theorem~\ref{thm:main} below that under mild technical conditions, if the $G$-varieties 
$X\times W$ and $X\times\field^{\dim(W)}$ are birationally isomorphic over $X$, then there is a $G$-equivariant 
birational isomorphism over $X$ between them which is {\it linear on} $W$, i.e. for a general $x\in X$, the restriction 
$\{x\}\times W\cong W\to \field^d\cong\{x\}\times \field^d$ is $\field$-linear. 

The $G$-variety $X$ is {\it generically free} if 
the stabilizer $G_x$ of a general $x\in X$ is trivial.  
The action of $G$ on $X$ is {\it generically separable} if the orbit morphism $G_0\to G_0x$, $g\mapsto gx$ is 
separable for a general $x\in X$ (this holds automatically when $\charac(\field)=0$).  
Note that if $X$ is a generically free $G$-variety, then generic separability is equivalent to the following: 
the map $G\times X\to X\times X$, $(g,x)\mapsto (x,gx)$ is birational between $G\times X$ and the graph 
$\{(x,gx)\mid g\in G, x\in X\}$ of the action. We refer to Chapter AG in \cite{borel} for the definition and basic properties of separability. 

Developing the idea of Lemmas~\ref{lemma:indep-cov}, \ref{lemma:converse-indep-cov}, and a result of Reichstein \cite{reichstein},   
we shall prove the following: 

\begin{theorem}\label{thm:main}  
Let $X$ be an irreducible generically separable $G$-variety, and $W$ a $d$-dimensional $G$-module. 
Then the following are equivalent: 
\begin{itemize}
\item[(1)] The stabilizer $G_x$ acts trivially on $W$ for a general $x\in X$. 
\item[(2)] There exist $d$ generically independent rational covariants from $X$ to $W$. 
\item[(3)] The $G$-varieties 
$X\times W$ and  $X\times \field^d$ (where $G$ acts trivially on $\field^d$) are birationally isomorphic over $X$. 
\item[(4)] There exists a $G$-equivariant birational isomorphism 
$X\times W\rto X\times \field^d$ (where $G$ acts trivially on $\field^d$) of the form 
$(x,w)\mapsto (x,\Phi(x,w))$ which is linear on $W$ (i.e. the map 
$\Phi(x,-):W\to\field^d$ is $\field$-linear for a general $x\in X$). 
\item[(5)] The invariant field $\field(X\times W)^G$ is purely transcendental 
over $\field(X)^G$ of transcendence degree $d$, generated by elements of the form $\sum_{j=1}^d\Phi_{ij}\varepsilon_j$, $(i=1,\ldots,d)$, 
where $\Phi_{ij}\in\field(X)$ and 
$\varepsilon_1,\ldots,\varepsilon_d$ is a basis in the dual space $W^*$ of $W$. 
\end{itemize}
\end{theorem} 

\begin{remark}\label{remark:precise} {\rm 
The implication $(1)\implies (3)$ (or some version of it) is a fundamental tool in the study of rationality questions. It has gone into the literature as the "no-name lemma" after \cite{dolgachev}. The first published reference (in characteristic zero, for a generically free $G$-variety $X$) is \cite{bogomolov}, 
see \cite{colliot-sansuc} for a survey. A proof for arbitrary $G$ is given in the recent paper 
\cite{gille}. A version for arbitrary base fields (using the concept of scheme theoretically free actions) is given in 
\cite{reichstein-vistoli}. Although these references consider generically free $G$-varieties $X$, 
in a remark attributed to Kraft, it is mentioned in \cite{dolgachev} that to make the conclusion (3), it is sufficient to assume that 
$G_x$ acts trivially on $W$ for a general $x\in X$. 

As far as we know, the fact that conclusion (3) can be strengthened to (4) and (5) has not been emphasized before, except in the case of finite groups 
(see the proof of Proposition 1.1 in \cite{endo-miyata}). 
For finite groups, this is essentially Speiser's Lemma 
(asserting that a finite dimensional $K$-vector space endowed with a 
semi-linear action of a finite group is spanned by invariant elements; 
see for example Lemma 2.3.8 in \cite{szamuely}) 
applied when $K$ is a function field of a $G$-variety. So one may view the implications $(1)\implies (4)$ and $(1)\implies (5)$ as a generalization of Speiser's Lemma for linear algebraic groups. 

The main point of our note is bringing (2) into the picture, and pointing out its equivalence to (4) and (5). 
Having this in mind, 
the implication $(1)\implies (2)$ (due to Reichstein \cite{reichstein} when $X$ is generically free and $\charac(\field)=0$) 
appears to be a more fundamental principle than the no-name lemma: 
its statement and proof are rather natural, and yield the latter in the above explicit form 
as an immediate corollary.  
Furthermore, the focus on covariants is a good starting point for a constructive approach; we shall comment on this later in Section~\ref{sec:reductive}.  
}\end{remark}

\begin{proof} The implication $(4)\implies (3)$ is trivial. 

$(3)\implies (1)$: The $G$-equivariant birational isomorphism $X\times W\rto X\times \field^d$ restricts to 
a $G_x$-equivariant birational isomorphism $\Phi(x,-):W\rto\field^d$ for a general $x\in X$; since $G_x$ acts trivially on $\field^d$, 
it acts trivially on $W$. 

$(2)\implies (4)$:  
Suppose $F_1,\ldots,F_d\in\rcov_G(X,W)$ are generically independent.  
Denote by $U$ the subset in $X$ where $F_1(x),\ldots,F_d(x)$ are all defined, it is a $G$-stable dense open subset in $X$. Apply Lemma~\ref{lemma:indep-cov} for the $G$-variety $U$. 

$(4)\implies (5)$ is straightforward: Suppose $(x,w)\mapsto (x,\Phi(x,w))$ is a $G$-equivariant 
birational isomorphism $X\times W\to X\times\field^d$, which is linear on $W$. 
Then $\Phi(x,w)=(\Phi_1(x,w),\ldots,\Phi_d(x,w))$ with 
$\Phi_i=\sum_{j=1}^d\Phi_{ij}\varepsilon_j\in\field(X\times W)^G$, where $\Phi_{ij}\in\field(X)$. 
Moreover, since $\field(X\times \field^d)^G$ is obviously generated over $\field(X)^G$ by the coordinate functions on $\field^d$, we get that $\field(X\times W)^G$ is generated over $\field(X)^G$ by the $d$ algebraically independent elements $\Phi_i$. 

$(5)\implies (2)$: Suppose that $\field(X\times W)^G$ is purely transcendental over $\field(X)^G$ generated by $\Phi_i=\sum_{j=1}^d\Phi_{ij}\varepsilon_j$, 
$i=1,\ldots,d$, where $\varepsilon_1,\ldots,\varepsilon_d$ is a basis of $W^*$. 
Write $g_W$ for the matrix of $g\in G$ acting on $W$ with respect to the basis 
$e_1,\ldots,e_d$ dual to $\varepsilon_i$. The $G$-invariance of the $\Phi_i$ can be expressed as the matrix equality 
\begin{equation}\label{eq:matrixeq} 
(g\Phi_{ij})_{d\times d}=(\Phi_{ij})_{d\times d}\cdot g_W. 
\end{equation} 
We claim that $\det(\Phi_{ij})_{d\times d}\neq 0\in\field(X)$. 
Indeed, assume on the contrary that the rows of $(\Phi_{ij})_{d\times d}$ are linearly dependent over $\field(X)$. After a possible reordering of the 
$\Phi_i$, we may assume that the first $r$ rows are linearly independent, whereas 
\begin{equation}\label{eq:lindep} 
(\Phi_{r+1,1},\ldots,\Phi_{r+1,d})=\sum_{i=1}^rf_i\cdot (\Phi_{i1},\ldots,\Phi_{id})
\end{equation} 
with $f_i\in\field(X)$. Next we show that all the $f_i$ are $G$-invariant: 
apply $g\in G$ to (\ref{eq:lindep}); Taking into account (\ref{eq:matrixeq}) and multiplying both sides of the resulting vector equality by $g_W^{-1}$ 
we obtain 
\begin{equation}\label{eq:glindep} 
(\Phi_{r+1,1},\ldots,\Phi_{r+1,d})=\sum_{i=1}^r(gf_j)
(\Phi_{i1},\ldots,\Phi_{id}). 
\end{equation} 
Take the difference of (\ref{eq:lindep}) and (\ref{eq:glindep}): 
\[(0,\ldots,0)=\sum_{i=1}^r(f_j-gf_j)
(\Phi_{i1},\ldots,\Phi_{id}). 
\]
Since the first $r$ rows of the matrix $(\Phi_{ij})_{d\times d}$ are linearly independent over $\field(X)$, it follows that $f_i=gf_i$ for $i=1,\ldots,r$. 
This holds for all $g\in G$, hence $f_1,\ldots,f_r\in\field(X)^G$. 
Consequently, we have 
\[\Phi_{r+1}=\sum_{i=1}^rf_i\Phi_i\in\field(X\times W)^G,\] 
contradicting the assumption that $\Phi_1,\ldots,\Phi_d$ are algebraically independent over $\field(X)^G$. 

Thus we proved that the matrix
$(\Phi_{ij})\in\field(X)^{d\times d}$ is invertible; denote by $F_{ij}\in\field(X)$ the $(i,j)$-entry of its inverse. 
Then formula (\ref{eq:matrixeq}) shows that 
$F_j:x\mapsto\sum_{i=1}^dF_{ij}(x)e_j$ 
is a covariant for $j=1,\ldots,d$. 
Moreover, these covariants are generically independent, since $\det(F_{ij})_{d\times d}=\det(\Phi_{ij})_{d\times d}^{-1}\neq 0$.  

$(1)\implies (2)$: Our first step is to reduce to the case when the action of $G$ on $X$ is generically free, and $W$ is a faithful $G$-module. 
Denote by $N$ the kernel of the action of $G$ on $W$. Then $N$ is a closed normal subgroup of $G$, and the stabilizer $G_x$ is contained in $N$ for a general $x\in X$, and $G/N$ acts faithfully on $W$. 
In fact we may assume $G_x\leq N$ for all $x$ by omitting a proper closed $G$-stable subset of $X$. 
Let $\pi:X\to X/N$ be a rational quotient; i.e., 
$X/N$ is a model (defined up to birational isomorphism) of $\field(X)^N$, and $\pi$ the dominant rational map corresponding to the field inclusion $\field(X)^N\to\field(X)$. There is a unique rational $G$-action (factoring through $G/N$) on $X/N$ such that 
$\pi$ is $G$-equivariant (see for example Theorem 5 in \cite{rosenlicht}). By a theorem of Weil \cite{weil} (see \cite{rosenlicht} for the case when $G$ is not connected) we may assume that $G/N$ acts morphically (not just rationally) on $X/N$. 

Now observe that the action of $G/N$ on $X/N$ is generically free. Indeed, let $U$ be an $N$-stable dense open subset of $X$ such that $\pi\vert_U:U\to\pi(U)=U/N$ is a geometric quotient morphism (i.e. it is an open morphism whose fibers are $N$-orbits); such an $U$ exists by Rosenlicht's Theorem \cite{rosenlicht}. 
It is easy to see that $\pi\vert_U$ extends to a $G$-equivariant morphism $\pi:\cup_{g\in G}gU\to\cup_{g\in G}g\pi(U)$, which is a geometric quotient with respect to the action of $N$. 
In other words, we may assume that $U$ is $G$-stable, hence passing from $X$ to a $G$-stable dense open subset if necessary, we may assume that the $G$-equivariant morphism $\pi:X\to X/N$ is a geometric quotient with respect to the action of $N$. 
Suppose $g\in G_{\pi(x)}$ for some $x\in X$. Then $gx=nx$ for some $n\in N$, hence $g^{-1}n\in G_x\leq N$, implying 
that $g\in N$. So for all $y\in X/N$ we have $G_y=N$. 

By Lemma~\ref{lemma:quotientsep} below, the action of $G/N$ on $X/N$ is generically separable.  

Note that if $F_1,\ldots,F_d\in\rcov_{G/N}(X/N,W)$ are generically independent, then 
$F_1\circ\pi,\ldots,F_d\circ\pi$ are generically independent rational covariants from $X$ to $W$. 
Therefore it is sufficient to deal with the case 
when $G$ acts generically freely and separably on $X$, and faithfully on $W$. From now on we assume that this is the case. In characteristic zero, the existence of $d$ generically independent covariants is due to Reichstein, see Lemma 7.4 in \cite{reichstein}. 
We follow the main idea of his proof, and extend it to positive characteristic.  
Identify $G$ with a subgroup of $H:=GL(W)$ in the obvious way. 
Then $H\times G$ acts on the product 
$H\times X$ via $(h,g)\cdot (h',x)=(hh'g^{-1},gx)$. Following Definition 2.12 in \cite{reichstein}, define 
$H\star_G X$ as a rational quotient of the $H\times G$-variety $H\times X$ 
with respect to the closed normal subgroup $G$ of $H\times G$; it is defined up to birational isomorphism, and is chosen to be an 
$H=(H\times G)/G$-variety; this is possible by the results from \cite{weil} and 
\cite{rosenlicht} mentioned above. 
(When $X$ is quasi-projective, $H\star_G X$ coincides with the corresponding homogeneous fiber product.)  We shall use the notation 
$[h,x]\in H\star_GX$ for $\pi(h,x)$.

Now $Y:=H\star_GX$ is a generically free $H$-variety, and $\iota: X\to Y$, $x\mapsto[1_H,x]$ is a $G$-equivariant embedding defined on a dense open subset of $X$. Moreover, $Y$ is a generically separable $H$-variety by Lemma~\ref{lemma:homfibprod}. 

Since $H$ is special (cf. Section 2.6 in  \cite{popov-vinberg}), 
there exists a rational section $\sigma:Y/H\to Y$ (i.e. setting 
$S:=\sigma(Y/H)$, we have that 
$H S$ is dense in $Y$, and $S$ intersects each $H$-orbit in at most one point). 
View $H\times S$ as an $H$-variety via $h\cdot(h's)=(hh',s)$. 
Then $H\times S\to Y$, $(h,s)\mapsto hs$ is an $H$-equivariant birational isomorphism, see for example Section 2.5 in \cite{popov-vinberg}, and  Lemma~\ref{lemma:section} for the case when $\charac(\field)>0$; denote by $\alpha$ its inverse. Take a point $x\in X$ such that 
$\iota(x)$ is in the domain of definition of $\alpha$ (there is such an $x$ since each $H$-orbit in $Y$ intersects $\iota(X)$, 
and the domain of definition of an $H$-equivariant rational map is $H$-stable). 
Denote by 
$\pr_1:H\times S\to H$ the projection morphism. Then $\pr_1\circ\alpha:Y\rto H$ is an $H$-equivariant rational map, and composing it with $\iota$ we 
get the composition 
\[X\stackrel{\iota}\rto Y\stackrel{\alpha}\rto H\times S\stackrel{\pr_1}\to H,\]   
which is a $G$-equivariant rational map $\varphi:X\rto H$ (defined at  $x$). 
For an arbitrary $w\in W$, we have the orbit morphism   
$\mu_{\varphi(x)^{-1}w}:H\to W$, $h\mapsto h\varphi(x)^{-1}w$. 
Clearly, $\eta_w=\mu_{\varphi(x)^{-1}w}\circ\varphi:X\to W$ is a $G$-equivariant rational map with $\eta_w(x)=w$. Taking a basis $e_1,\ldots,e_d$ of $W$, then 
$\eta_{e_i}\in\rcov_G(X,W)$, $i=1,\ldots,d$, are linearly independent at $x\in X$. 
\end{proof}

The following four Lemmas were used in the preceeding proof. These technical statements relating mainly separability must be well known. 
We include their proof, since we did not find a convenient reference. 

\begin{lemma}\label{lemma:quotientsep} 
Let $X$ be an irreducible $G$-variety, $N$ a closed normal subgroup of $G$ such that $G_x\leq N$ for all $x\in X$, and let 
$\pi:X\to Y$ be a morphism of $G$-varieties, which is a rational  quotient with respect to the action of $N$.  
If the $G$-variety $X$ is generically separable, then the $G/N$-variety $Y$ is generically separable. 
\end{lemma} 

\begin{proof} Write $\overline{G}:=G/N$. We may assume that the fibers of $\pi$ are $N$-orbits (see the proof of $(1)\implies (2)$ in Theorem~\ref{thm:main}), and so $\overline{G}$ acts freely on $Y$. The field 
$\field(X)$ is a separable extension of $\field(X)^N$ (see for example \cite{borel}), so 
the differential $\dif_x\pi:T_xX\to T_{\pi(x)}Y$ is surjective for a general $x\in X$, hence 
\[\dim(\ker(\dif_x\pi))=\dim(X)-\dim(Y)=\dim(Nx).\] 
On the other hand,  
$\ker(\dif_x\pi)\supseteq\ker(\dif_x\pi\vert_{T_xGx})\supseteq T_x(Nx)$, 
and $\dim(T_xNx)=\dim(Nx)$. It follows that 
\[\ker(\dif_x(\pi))=T_x(Nx)=\ker(\dif_x\pi\vert_{T_xGx})=\ker(\dif_x(\pi\vert_{Gx})).\]  
Hence (for a general $x\in X$) we have \[\dim(\imm(\dif_x\pi\vert_{Gx}))=\dim(Gx)-\dim(Nx)=\dim(\overline{G}\pi(x))\]  (in the last equality we use that the stabilizer of $\pi(x)$ in $G$ is $N$). 
So $\dif_x(\pi\vert_{Gx}):T_xGx\to T_{\pi(x)}\overline{G}\pi(x)$ 
is surjective. The differential of the morphism $G\to Gx$, $g\mapsto gx$ is also surjective (by the assumption on generic separability of the action of $G$ on $X$), hence 
the differential of $G\to \overline{G}\pi(x)$, $g\mapsto g\pi(x)$ is surjective at $1\in G$. 
This morphism is the composition of the natural surjection $G\to \overline{G}$ and the morphism $\overline{G}\to \overline{G}\pi(x)$, $g\mapsto g\pi(x)$, implying that the differential at $1\in \overline{G}$ of the latter is surjective. This holds for a general point $\pi(x)\in Y$, thus the action of $\overline{G}$ on $Y$ is indeed generically separable. 
\end{proof} 

\begin{lemma}\label{lemma:homfibprod} 
Let $X$ be a generically free and separable irreducible $H$-variety, where $H$ is a closed subgroup of the linear algebraic group 
$G$. Then 
$G\star_HX$ is a generically separable $G$-variety. 
\end{lemma} 

\begin{proof}  Consider $G\times X$ together with the 
$G\times H$-action $(g,h)\cdot(g',x)=(gg'h^{-1},hx)$. The tangent space $T_{(g,x)}(G\times X)$ decomposes as 
$T_g(G)\oplus T_x(X)$. Denote by $\mu_{(g,x)}$ the orbit morphism  $G\times H\to G\times X$ at $(g,x)$. The differential of its restriction to $G$ maps $T_1G$ isomorphically onto $T_gG\subset T_{(g,x)}(G\times X)$, and the projection to $T_xX$ of the image of 
$\dif_{1_H}\mu\vert_H$ is $T_xX$ for a general $x\in X$ (since $X$ is a generically separable $H$-variety). 
This shows that $G\times X$ is a generically separable $G\times H$-variety. It is obviously free. 
By definition, the $G=(G\times H)/H$-variety $G\star_HX$ is a rational quotient of $G\times X$ with respect to the action of the normal subgroup $H$ 
in $G\times H$. So generic separability of the $G$-variety  $G\star_HX$ follows from Lemma~\ref{lemma:quotientsep}. 
\end{proof}

\begin{lemma}\label{lemma:section} Let $X$ be a generically free and separable irreducible $G$-variety, where 
$G$ is the general linear group $GL_n(\field)$. Then there is a $G$-equivariant birational isomorphism 
between $X$ and $G\times S$ for some variety $S$ 
(where $G$ acts on $G\times S$ via $g\cdot(g',s)=(gg',s)$). 
\end{lemma} 

\begin{proof} The Galois cohomology $H^1(\Gamma,GL_n(K))$ is trivial for any finite group $\Gamma$ of automorphisms of a 
field $K$, see for example the section on Galois descent in \cite{szamuely}. When $\charac(\field)=0$, 
it is explained in Section 2.5 in \cite{popov-vinberg} how to derive our statement from this fact. 
We adjust this argument to positive characteristic. 

Replacing $X$ by a $G$-stable dense open subset, we may assume that a geometric quotient 
$\pi:X\to X/G=:Y$ exists. 
Identify $\field(Y)$ with the subfield $\pi^*(\field(Y))$ in $\field(X)$. 
Let $U$ be a dense affine open subset of $X$. 
Denote by $A$ the $\field(Y)$ subalgebra in $\field(X)$ generated by $\field[U]$. 
A $\field(Y)$-algebra homomorphism $A\to\field(Y)$ restricts to a $\field$-algebra homomorphism of $\field[U]$ 
onto a finitely generated $\field$-subalgebra of $\field(Y)$; this $\field$-algebra homomorphism is the comorphism of a morphism 
$\sigma:V\to U$ (where $V$ is some dense affine open subset of $Y$). Clearly, $\pi\circ\sigma=\iid_{V}$. 
So to a $\field(Y)$-point of $\spec(A)$ there corresponds a rational section 
$\sigma:Y\rto U$ for $\pi$. 

Recall that a quasisection is a rational map $\sigma:Z\to X$ such that $\pi\circ\sigma:Z\rto Y$ is dominant and finite, 
(i.e. $\field(Z)$ is a finite extension of $\field(Y)$). 
Fix a finite extension $\field(Z)$ of $\field(Y)$ (i.e. a dominant finite rational map $\alpha: Z\rto Y$. In the same way as above, to a 
$\field(Z)$-point of $\spec(A)$ corresponds a quasisection $\sigma:Z\rto U$ with 
$\pi\circ\sigma=\alpha$. 

Now $\spec(A)$ has a point in the separable closure of $\field(Y)$ (see for example Chapter AG in \cite{borel} or the Appendix of \cite{szamuely}), 
hence $\spec(A)$ has a point already in a finite Galois extension $\field(Z)$ of $\field(Y)$ (write $\alpha:Z\to Y$ for the corresponding finite, dominant, rational map; we may assume that $\alpha$ is an everywhere defined morphism). 
Denote by $\sigma:Z\to X$ the corresponding quasisection (we may assume that $\dom(\sigma)=Z$). 
Consider the fibre product $\tilde X:=X\times_YZ$, the pull back $\tilde\pi:\tilde X\to Z$ of $\pi$. 
and $\tilde\sigma:Z\to \tilde X$ corresponding to $\sigma$. Then $\tilde\pi\circ\tilde\sigma=\iid_Z$, so $\tilde\sigma$ is a section for $\tilde\pi$.  

$\tilde X$ is naturally a $G$-variety, it is obviously generically free and separable. Therefore by Lemma~\ref{lemma:biratproduct} below, there is a $G$-equivariant birational isomorphism $\beta:\tilde X\to G\times Z$ ($G$ acts trivially on $Z$ and by left multiplication on itself). 
Now let $\gamma$ be an element of the Galois group $\Gamma$ of the Galois extension $\field(Z)\supset\field(Y)$. Write $\gamma^*$ for the birational isomorphism $Z\rto Z$ corresponding to $\gamma$. 
We have $\alpha\circ\gamma^*=\alpha$. Then $\tilde\sigma^{\gamma}:=\tilde\sigma\circ\gamma^*:Z\to\tilde X$ is another rational section of $\tilde\pi$. 
Moreover, for a general $z\in Z$, there is a unique $\sigma_{\gamma}(z)\in G$ with $\tilde\sigma^{\gamma}(z)=\sigma_{\gamma}(z)\tilde\sigma(z)$. 
Here $\sigma_{\gamma}:Z\to G$ is rational, since the mapping 
$\beta\circ\tilde\sigma\circ\gamma^*$ is rational, and it agrees with $z\mapsto(\sigma_{\gamma}(z),z)$. 
So one may view $\sigma_{\gamma}$ as an element of $GL_n(\field(Z))$. 
One checks easily that the map $\Gamma\to GL_n(\field(Z))$, $\gamma\mapsto\sigma_{\gamma}$ is a 1-cocycle, i.e. 
$\sigma_{\gamma\rho}=(\sigma_{\gamma})^{\rho}\sigma_{\rho}$ (we extend the action of $\Gamma$ to $GL_n(\field(Z))$). 
Triviality of $H^1(\Gamma,GL_n(\field(Z))$ means that there exists a $g\in GL_n(\field(Z))$ such that denoting by 
$\tau$ the quasisection $z\mapsto g(z)\sigma(z)$, we have that $\tau^{\gamma}=\tau$ for all $\gamma\in\Gamma$. 
This means that the quasisection $\tau$ corresponds to a point of $\spec(A)$ in $\field(Z)^{\Gamma}=\field(Y)$, i.e., 
$\tau$ is a rational section of $\pi$.  
\end{proof} 

\begin{lemma}\label{lemma:biratproduct} Let $X$ be a a generically free and separable $G$-variety, 
$\pi:X\to Y$ a morphism whose fibres are $G$-orbits, and $\sigma:Y\to X$ is a morphism with 
$\pi\circ\sigma=\iid_Y$. Then $G\times\sigma(S)\to X$, $(g,s)\mapsto gs$ is a birational isomorphism. 
\end{lemma} 

\begin{proof} Set $S:=\sigma(Y)$, and $\mu:G\times S\to X$, $(g,s)\mapsto gs$. 
For any $s\in S$, the differential $\dif_s\pi\vert_S$ maps $T_sS$ isomorphically onto 
$T_{\pi(s)}Y$ (since $\pi\vert_S:S\to Y$ is an isomorphism). Moreover, $\ker(\dif_s\pi)\supseteq T_sGs$. 
For a general $s\in S$ we have $\dim(X)=\dim(Gs)+\dim(T_sS)$, implying that 
$T_sX=T_sGs\oplus T_sS$. For such an $s$, 
consider $\dif_{(1_G,s)}\mu:T_{(1_G,s)}G\times S\to T_sX$. 
By assumption, the orbit morphism $\mu_s:G\to Gs$, $g\mapsto gs$ is separable for a general $s\in S$ (indeed, 
$X=GS$, and $\mu_s$, is separable if and only if $\mu_{gs}$ is separable for some $g\in G$). 
It follows that the restriction of $\dif_{(1_G,s)}\mu$ to $T_{1_G}G$ (we identify $T_{(1_G,s)}G\times S$ with 
$T_{1_G}G\oplus T_sS$) maps surjectively onto $T_sGs$, whereas the image of the restriction of 
$\dif_{(1_G,s)}\mu\vert_{1_G\times S}$ to $T_{(1_G,s)}(1_G\times S)$ is $T_sS$. This shows that $\dif_{(1_G,s)}\mu$ is surjective onto $T_sX$. 
This holds for a general $s\in S$, showing that $\mu$ is a separable morphism. 
The generic freeness of the action implies that $\mu$ is bijective. 
It follows that $\mu$ is birational, see eg. Chapter AG in \cite{borel}.  
\end{proof}


\section{Linearly reductive groups} \label{sec:reductive} 

For linearly reductive groups, the results of Section~\ref{sec:mainlemma} 
have a counterpart using integral covariants, that fits better into the traditions of invariant theory. 
Let us recall the following corollary of Theorem 1 from \cite{panyushev} (some special cases appear in \cite{lehrer-springer} and Proposition 4.2.5 in \cite{derksen-kemper}): 

\begin{proposition}\label{prop:panyushev} 
Let $G$ be a linearly reductive group, 
$X$ an affine $G$-variety, and $W$ a $G$-module. 
If for some $x\in X$ having closed orbit the stabilizer $G_x$ acts trivially on $W$, then there exist 
$d=\dim_{\field}(W)$ covariants $F_1,\ldots,F_d\in\cov_G(X,W)$ 
such that $F_1(x),\ldots,F_d(x)$ are linearly independent over $\field$. 
\end{proposition} 

\begin{proof} For an arbitrary $w\in W$ the map 
$F_w:gx\mapsto gw$ is obviously a covariant from the orbit $Gx\cong G/G_x$ to $W$. 
The coordinate ring of $Gx$ can be viewed as a $G$-module direct summand of 
$\field[X]$ (the linear reductivity of $G$ is used here), therefore $F_w$ can be lifted to a covariant $H_w:X\to W$, with $H_w(x)=w$. Since this holds for all $w\in W$, the claim follows. 
\end{proof} 

\begin{remark}\label{remark:stabilizer} 
{\rm The assumption that "$G_x$ acts trivially on $W$" is necessary in Proposition~\ref{prop:panyushev}, since the covariance of the $F_i$ implies  that $F_1(x),\ldots,F_d(x)$ are all fixed by $G_x$, and these vectors constitute a basis of $W$. }
\end{remark}

\begin{remark}\label{remark:panyushev} 
{\rm As it is noted in \cite{panyushev}, instead of the closedness of the orbit $Gx$, it is sufficient to assume that the closure of $Gx$ is normal and the complement of $Gx$ in its closure has codimension $\geq 2$ } 
\end{remark}

\begin{corollary}\label{cor:reductive} 
Let $G$ be a linearly reductive group, 
$X$ an irreducible affine $G$-variety, and $W$ a $G$-module. 
If for some $x\in X$ having closed orbit the stabilizer $G_x$ acts trivially on $W$, then there exists a relative invariant $f\in\field[X]$ such that $f(x)\neq 0$ and the $G$-varieties  $X_f\times W\to X_f\times\field^{\dim(W)}$ 
($G$ acting trivially on $\field^{\dim(W)}$) 
are isomorphic over $X$. 
\end{corollary} 

\begin{proof} Immediate from Proposition~\ref{prop:panyushev} and Lemma~\ref{lemma:indep-cov}. 
\end{proof} 
\begin{remark}\label{remark:constructive} {\rm 
Corollary~\ref{cor:reductive} can be viewed as a constructive version of the 
no-name lemma. Indeed, if $X$ is a $G$-module, then an algorithm to compute $\field[X]^G$-module generators of $\cov_G(X,W)$ is explained on page 157 of 
\cite{derksen-kemper}. Then it is an easy matter to select 
$d=\dim_{\field}(W)$ of the generators that are generically independent (or to show that there are no $d$ generically independent covariants), and get $f$ and the isomorphism $X_f\times W\to X_f\times\field^d$ explicitly as in the proof of 
Lemma~\ref{lemma:indep-cov}. }
\end{remark} 

\begin{corollary}\label{cor:2reductive} 
Let $G$ be a linearly reductive group, 
$X$ an irreducible affine $G$-variety, and $W$ a $G$-module. The following are equivalent: 
\begin{itemize} 
\item[(1)] There is a non-zero relative invariant $f\in\field[X]$ such that $G_x$ acts trivially on $W$ for all 
$x\in X_f$. 
\item[(2)] There is a non-zero relative invariant $f\in\field[X]$ and an $x\in X_f$ whose $G$-orbit is closed in $X_f$, and 
$G_x$ acts trivially on $W$. 
\item[(3)] There is a non-zero relative invariant $f\in\field[X]$ such that there are 
$d:=\dim(W)$ covariants 
$F_i:X_f\to W$ $(i=1,\ldots,d)$ with $F_1(x),\ldots,F_d(x)$ linearly independent for all $x\in X_f$. 
\item[(4)] There is a non-zero relative invariant $f\in\field[X]$ such that the $G$-varieties $X_f\times W$ and $X_f\times \field^{\dim(W)}$ 
are isomorphic. 
\end{itemize} 
\end{corollary} 

\begin{proof} $(1)\implies (2)$: The $G$-variety $X_f$ contains a point $x$ whose $G$-orbit is closed. 

$(2)\implies (3)$: Apply Proposition~\ref{prop:panyushev} and the proof of Lemma~\ref{lemma:indep-cov}.  

$(3)\implies (4)$: Apply Lemma~\ref{lemma:indep-cov}. 

$(4)\implies (1)$: The $G$-equivariant isomorphism $X_f\times W\to X_f\times \field^d$ restricts for any $x\in X_f$ to a $G_x$-equivariant isomorphism 
$\{x\}\times W\to \{x\}\times \field^d$. Now $G_x$ acts trivially on the latter, hence $G_x$ acts trivially on $W$.  
\end{proof} 

\begin{remark}\label{remark:shmelkin} {\rm In general in an affine $G$-variety $X$ there are more closed $G'$-orbits than $G$-orbits, 
where $G'$ is the commutator subgroup of $G$. 
When $\charac(\field)=0$ and $G$ is connected reductive, by Theorem 1 in \cite{shmelkin}, $G'x$ is closed in $X$ if and only if there is a relative invariant  $b$ on $X$ such that $b(x)\neq 0$ and the $G$-orbit of $x$ is closed in the affine open set $X_b$. 
Therefore in this case (1), (2), and (3) in Corollary~\ref{cor:2reductive} are equivalent also to the following:  

(4)  There is a point $x\in X$ whose $G'$-orbit is closed, 
and the stabilizer $G_x$ of $x$ in $G$ acts trivially on $W$. 
} \end{remark}

The considerations above can be applied to generically free actions of linearly reductive groups on factorial affine varieties: 
(An affine variety $X$ is called {\it factorial} if $\field[X]$ is a unique factorization domain.) 

\begin{theorem}\label{thm:almost-free} 
Suppose $\charac(\field)=0$, $G$ is reductive, acting generically freely on the factorial affine variety $X$.  Then for any $G$-module $W$, there is a non-zero relative invariant $f$ on $X$ such that the $G$-varieties $X_f\times W$ and 
$X_f\times\field^d$ are isomorphic over $X_f$. 
\end{theorem} 

\begin{proof} By a result of M. Van den Bergh (see Theorem 2.18 in \cite{colliot-sansuc}), there is a non-empty affine open $G$-stable subset $U$ in $X$ such that the generic $G$-orbit is closed in $U$. 
Clearly, $U$ is of the form $X_b$ for some non-zero relative invariant $b$ 
(see Theorem 3.1 in \cite{popov-vinberg}).  
Now apply Corollary~\ref{cor:reductive} for the affine $G$-variety $U$ and $W$. 
\end{proof} 

\begin{remark} {\rm We conclude this section by drawing attention to another application of covariants, which has some common flavour with ours. In 4.2.1 in \cite{derksen-kemper}, 
covariants are applied to refine a corollary of Luna's slice theorem (in a much more special situation). } 
\end{remark}


\section{On a theorem of Molien}\label{sec:molien} 

If $G$ is a finite group acting faithfully on $X$, and $|G|$ is invertible in $\field$, then the assumptions of 
Proposition~\ref{prop:panyushev} hold for all $G$-modules $W$, hence there are 
$\dim(W)$ generically independent covariants from $X$ to $W$. 
(The weaker statement that for any finite subgroup of $GL(V)$, the polynomial ring $\field[V]$ contains all irreducible $G$-modules $W$ as a summand is usually refered to as Molien's Theorem.) 
We take a digression now and extend this result to the modular case: 

\begin{proposition}\label{prop:modular} 
Let $G$ be a finite group acting faithfully on a factorial affine variety $X$ (over a base field $\field$ of arbitrary 
characteristic). Then for any $G$-module $W$, there are $\dim(W)$ generically independent covariants from 
$X$ to $W$. 
\end{proposition}

\begin{proof} By Theorem~\ref{thm:main} $(1)\implies (2)$ (or by Speiser's Lemma and the implication $(5)\implies (2)$ in Theorem~\ref{thm:main}) we know that there are $d:=\dim(W)$ 
generically independent rational covariants $F_i:X\rto W$, $i=1,\ldots,d$. 
Take a non-zero $h\in \field[X]$ such that the $F_i$ are all defined on $X_h$. 
Set $f:=\prod_{g\in G}g\cdot h$. Then $f\in\field[X]^G$ is non-zero, and the coordinate functions $F_{ij}$ of the $F_i$ belong to 
$\field[X_f]$ for all $i,j=1,\ldots,d$. Clearly there exists an exponent $n$ such that $f^nF_{ij}\in\field[X]$ for all 
$i,j$. Then  $f^nF_i\in \cov_G(X,W)$, $i=1,\ldots,d$ are generically independent integral covariants from $X$ to $W$.  
\end{proof}


\section{Examples}\label{sec:examples} 

In this section we collect a couple of situations where there are some "obvious" covariants, and 
our approach yields explicit birational isomorphisms and explicit generators of fields of invariants.  
The following trivial lemma makes possible a successive application of Lemma~\ref{lemma:indep-cov} in many situations. 

\begin{lemma}\label{lemma:projection} 
Let $X$ be an irreducible $G$-variety and $W$ a $G$-module, and assume that $F_1,\ldots,F_d\in \cov_G(X,W)$ are generically independent covariants. 
For an arbitrary irreducible $G$-variety $Y$ denote $\pi_X:X\times Y\to X$  
the projection morphism. 
Then 
$F_1\circ\pi_X,\ldots,F_d\circ\pi_X$ are generically independent covariants 
from $X\times Y$ to $W$. 
\end{lemma} 

\begin{example}\label{example:algebras} {\rm Let $R$ be a finite dimensional associative (or Lie or Jordan) algebra over an infinite  field $\field$ and $G$ a group of $\field$-algebra automorphisms of $R$. 
Suppose that the algebra $R$ is generated by $n$ elements. 
Take the free associative $\field$-algebra on $n$ generators $x_1,\ldots,x_n$. 
Fix a generating system $a_1,\ldots,a_n$ of $R$. Choose $d:=\dim_{\field}(R)$ monomials $M_i=M_i(x_1,\ldots,x_n)$ ($i=1,\ldots,d$) in the free algebra such that $M_i(a_1,\ldots,a_n)$, $i=1,\ldots,d$, is a basis of $R$. 
Identify $M_i$ with the function $R^n\to R$, 
$(r_1,\ldots,r_n)\mapsto M_i(r_1,\ldots,r_n)$. 
Then $M_i\in\cov_G(R^n,R)$ for $i=1,\ldots,d$, and these covariants are generically independent by construction. 
In particular, by Lemma~\ref{lemma:indep-cov}, 
for any integer $m\geq n$, the field 
$\field(R^m)^G$ of $G$-invariant rational functions on $R^m=R\times\cdots\times R$ is purely transcendental over $\field(R^n)^G$ of transcendence degree $(m-n)\dim_{\field}(R)$. 
Since there is a general interest in classifying $m$-tuples of elements in an algebra $R$ up to a group of automorphisms of 
$R$, this example seems worthwile to be mentioned. 

A notable special case is when 
$R=M(n,\field)$ is the algebra of $n\times n$ matrices over an arbitrary infinite base field, and  $G:=GL_n(\field)$ acts on $R$ by 
conjugation. 
Consider the diagonal action of $G$ on $R^m:=R\oplus\cdots \oplus R$, the space of $m$-tuples of $n\times n$ matrices. 
It is a long standing open problem whether the field $\field(V^m)^G$ 
is rational (a purely transcendental extension of $\field$), see \cite{colliot-sansuc} for references. 
The answer is known to be positive in a few special cases only. 
It was an essential observation in \cite{procesi} that for $m\geq 2$, the  field $\field(R^m)^G$ is rational over $\field(V^2)^G$. 
This has an easy explanation using the following covariants in $\cov_G(R^2,R)$: 
$(A,B)\mapsto A^iB^j$, where $A,B\in M(n,\field)$, and $i,j=0,\ldots,n-1$. 
It is easy to see that these $n^2$ covariants are generically independent. 
(Indeed, substitute for $A$ a diagonal matrix with $n$ distinct diagonal entries, and for $B$ the permutation matrix corresponding to an $n$-cycle.) 

Moreover, one gets automatically the 
rationality of $\field(R^m)^H$ over $\field(R^2)^H$ for any subgroup $H$ of 
$G$. 
The case when $H$ is the  
orthogonal group ${\mathrm{O}}(n,\field)$ 
or the symplectic group ${\mathrm{Sp}}(n,\field)$ (when $n$ is even). 
is studied in \cite{saltman} over an infinite field of characteristic different from $2$,  
and the rationality of $\field(R^m)^H$ over $\field(R^2)^H$ is proved there in a bit more elaborate way.  }
\end{example}

\begin{example}\label{example:weyl} {\rm 
Let $V$ be an $n$-dimensional $G$-module. 
Then the maps $V^n\to V$, 
$(v_1,\ldots,v_n)\mapsto v_i$, 
$i=1,\ldots,n$ are generically independent covariants from $V^n$ to $V$. 
Let $x_{ij}$ denote the coordinate function on $V^m$ mapping an $m$-tuple of vectors to the $i$th coordinate of the $j$th component, write 
${\underline{x}}_j$ for the column vector $[x_{1j},\ldots,x_{nj}]^T$, and consider the 
$n\times n$ matrix $D:=(x_{ij})_{i,j=1}^n$. 
For an arbitrary irreducible $G$-variety $X$ and $m\geq n$ 
we have by the proof of Lemma~\ref{lemma:indep-cov} and by Lemma~\ref{lemma:projection} that 
$\field(X\times V^m)^G$ is purely transcendental over $\field(X\times V^n)^G$ 
generated by the entries of $D^{-1}{\underline{x}}_i$, $i=n+1,\ldots,m$. 
More precisely, the localization 
$\field[X\times V^m]^G_{\det(D)}$ is a polynomial ring over 
$\field[X\times V^n]^G_{\det(D)}$ generated 
by the entries of $D^{-1}{\underline{x}}_i$, $i=n+1,\ldots,m$. 
This result is contained implicitly for example in the proof of Theorem 3 in \cite{grosshans:2007}. } 
\end{example} 

\begin{example} \label{example:permutation} {\rm 
Another case when one can easily write down covariants is that of a permutation action, i.e. let  the finite group $G$ 
act on $V:=\field^n$ by permutating the coordinates. 
Then 
\[\left[\begin{array}{c}x_1 \\\vdots \\x_n\end{array}\right] 
\mapsto \left[\begin{array}{c}x_1^i \\\vdots \\x_n^i\end{array}\right], 
\quad i=1,\ldots,n\] 
are generically independent covariants from $V\to V$, which can be used to 
construct an explicit birational isomorphism 
$V^m\to V\times \field^{(m-1)n}$. 
Just like in Example~\ref{example:weyl}, write $x_{ij}$ for the coordinate functions  
on $V^m$, and write 
$\Delta$ for the Vandermonde matrix $(x_{i}^j)_{i,j=1}^n$. 
Then the localization 
$\field[V^m]^G_{\det(\Delta)}$ is a polynomial ring over 
$\field[V]^G_{\det(\Delta)}$ generated by the entries of 
$\Delta^{-1}{\underline{x}}_j$, $j=2,\ldots,m$.}  
\end{example} 


\section{On generic independence of covariants}\label{sec:gen-indep} 

In this section we make some observations on generic independence of covariants. Let $X$ be an irreducible $G$-variety and 
$W$ a $G$-module. 
Denote by $\mor(X,W)$ (respectively $\rat(X,W)$) the $\field$-vector space of morphisms (respectively rational maps) $X\to W$ of  
algebraic varieties. $\rat(X,W)$ is naturally a $\dim_{\field}(W)$-dimensional vector space over $\field(X)$,  
and $\mor(X,W)$ is a free $\field[X]$-module of rank $\dim_{\field}(W)$. 
The group $G$ acts on $\rat(X,W)$ as follows: for $g\in G$, $F\in\mor(X,W)$, and $x\in X$ 
we have $(g F)(x)=g\cdot F(g^{-1}x)$. Clearly 
$\rcov_G(X,W)\subseteq\rat(X,W)$ is the subset of fixed points under this action. Moreover, $\mor(X,W)$ is a $G$-stable subset of $\rat(X,W)$, and $\cov_G(X,W)$ is the subset of fix points. 

By elementary linear algebra, $F_1,\ldots,F_e\in \rcov_G(X,W)$ are generically independent if and only if they are independent over $\field(X)$ in 
$\rat(X,W)$. One has therefore the following: 

\begin{proposition}\label{prop:indep}
If $F_1,\ldots,F_e\in\cov_G(X,W)$ are generically independent, then 
they are independent over the ring of invariants 
$\field[X]^G$ in the $\field[X]^G$-module $\cov_G(X,W)$. 
\end{proposition} 

The converse of Proposition~\ref{prop:indep} does not hold in general. 
For example, the group $G:=\field^{\times}$ acts on $X:=\field^2$ and $W:=\field$ by scalar multiplication. 
Then the covariants $(x,y)\mapsto x$ and $(x,y)\mapsto y$ from $X$ to $W$ are independent over 
$\field[X]^G=\field$, but they are not generically independent (since the number of generically independent covariants with values in $W$ is bounded by $\dim(W)$). 

However, working in the framework of rational covariants, Reichstein \cite{reichstein} proved the following: 

\begin{proposition} \label{prop:reichstein} {\rm (Lemma 7.4 (a) in \cite{reichstein})} 
The rational covariants $F_1,\ldots,F_e\in\rcov_G(X,W)$ are generically independent if and only if 
they are linearly independent over $\field(X)^G$ in $\rcov_G(X,W)$. 
\end{proposition} 

This yields the following partial converses for Proposition~\ref{prop:indep}: 

\begin{corollary}\label{cor:integralindep} 
If $\field(X)^G$ is the field of fractions of $\field[X]^G$ (this holds for example when $G$ is finite and $X$ is affine, or when 
$G=G'$, $X$ is factorial affine and $\field[X]$ does not contain non-scalar units). 
Then $F_1,\ldots,F_e\in\cov_G(X,W)$ are generically independent if and only if they are independent in the $\field[X]^G$-module  
$\cov_G(X,W)$. 
\end{corollary} 

\begin{corollary}\label{cor:relativeinvariant} 
Suppose $X$ is a factorial affine variety with $\field[X]^{\times}=\field^{\times}$, $F_1,\ldots,F_e\in\cov_G(X,W)$, 
such that $F_1(x),\ldots,F_e(x)$ are linearly dependent over $\field$ for all $x\in X$.   
Then there are relative invariants $h_1,\ldots,h_e\in\field[X]$ of the same weight such that 
$h_1F_1+\cdots +h_eF_e=0$ in $\mor(X,W)$. 
\end{corollary} 

\begin{proof} The assumptions on $X$ guarantee that any rational invariant in $\field(X)^G$ can be written as a fraction 
of relative invariants having the same weight. 
Now by Proposition~\ref{prop:reichstein} we know that there exist  $a_i\in\field(X)^G$, $i=1,\ldots,e$, such that 
$\sum_{i=1}^ea_iF_i=0$. Multiply this equality by $\prod_{i=1}^ec_i$, where $a_i=b_i/c_i$, and $b_i,c_i$ are relative invariants in 
$\field[X]$ of the same weight.  
\end{proof} 

The above two results maybe useful, since (especially for linearly reductive groups) 
much information on covariants is encoded in $\field[X]$ viewed as a 
$G-\field[X]^G$-module (see for example 3.13 in \cite{popov-vinberg}),  
and there are situations when information is available on the $\field[X]^G$-module structure of 
$\cov_G(X,W)$. 

Finallly we point out one more situation when the converse of Proposition~\ref{prop:indep} holds: 

\begin{proposition}\label{prop:pseudo-refl-indep} 
Assume that the linear algebraic group $G\subset\gl(V)$ is generated as an algebraic group by elements fixing pointwise a hyperplane in $V$, and let $W$ be a $G$-module. 
Some covariants $F_1,\ldots,F_e\in\cov_G(V,W)$ are generically independent if and only if they are independent 
in the $\field[V]^G$-module $\cov_G(V,W)$. 
\end{proposition} 

\begin{proof} The "only if" part follows from Proposition~\ref{prop:indep}. 
We prove the "if" part. Assume that $F_1,\ldots,F_e$ are covariants independent over $\field[V]^G$. 
As we noted at the beginning of this section, it is sufficient to show that they are independent over $\field[V]$ (and hence over $\field(V)$). Assume 
\begin{equation}\label{eq:0707270} 
h_1F_1+\cdots+h_eF_e=0\end{equation} 
for some 
$h_i\in\field[V]$, not all zero. 
Then for all $v\in V$ and $g\in G$ we have 
\[0=g\cdot \sum_{i=1}^eh(g^{-1}v)F_i(g^{-1}v)
=g\cdot \sum_{i=1}^eh(g^{-1}v)g^{-1}\cdot F_i(v)=\sum_{i=1}^eh_i(g^{-1}v)F_i(v),\] 
hence 
\begin{equation}\label{eq:6}\sum_{j=1}^e(g\cdot h_j)F_j=0\end{equation} 
holds for all $g\in G$. 

Choose the relation (\ref{eq:0707270}) so that the maximum degree of the $h_i$ on the left hand side 
is minimal possible.  
Now let $g$ be an element in $G$ fixing the hyperplane 
$\ker(\lll)$, where $\lll\in V^*$, the dual space of $V$. 
Subtract from (\ref{eq:0707270}) the equation 
(\ref{eq:6}) to get 
\[\sum_{j=1}^e(h_j-g\cdot h_j)F_j.\] 
The polynomial $h_j-g\cdot h_j$ vanishes on the hyperplane 
$\ker(\lll)$, hence $h_j-g\cdot h_j=\lll b_j$ for some 
$b_j\in\field[V]$. 
Dividing the above relation by $\lll$ 
we conclude that 
$\sum_{j=1}^eb_jF_j=0$. 
If some $b_j\neq 0$, then this is a nontrivial relation 
of strictly lower degree than (\ref{eq:0707270}), contradicting the minimality of the degree of the chosen relation. 
Therefore all $b_j=0$. 
This means that $h_j-g\cdot h_j=0$ for all elements $g\in G$ that fix some  hyperplane. 
Since $G$ contains a Zariski dense subgroup generated by such elements, it follows that 
$h_j$ is a $G$-invariant for $j=1,\ldots,e$. 
This contradicts the assumption that the $F_j$ are independent over 
$\field[V]^G$. 
\end{proof} 

\begin{example}\label{example:pseudoref} {\rm 
(i) A particularly nice situation when the conditions of the above proposition 
apply is when $G\subset GL(V)$ is a complex pseudo-reflection group 
(i.e. $\field={\mathbb{C}}$, the field of complex numbers, and $G$ is a finite group generated by elements fixing pointwise a hyperplane in the complex vector space $V$). Then $\field[V]^G$ is a polynomial ring and 
$\field[V]$ is a free $\field[V]^G$-module, see \cite{chevalley}. (We mention that the assumption on $G$ is used in a similar manner in the proof of Proposition~\ref{prop:pseudo-refl-indep} as in the proof of Lemma 1 in \cite{chevalley}.) 
Moreover, if $W$ is an irreducible $G$-module, then 
$\cov_G(V,W)$ is a free $\field[V]^G$-module of rank $\dim_{\field}(W)$ 
(see \cite{chevalley} and the discussion in 3.13 in \cite{popov-vinberg}), so one can find $\dim_{\field}(W)$ generically independent covariants by decomposing the coinvariant algebra (the quotient of $\field[V]$ modulo the ideal generated by the homogeneous $G$-invariants of positive degree). 

(ii) Over any base field, the general linear group $GL(V)$, the special linear group $SL(V)$, 
the full orthogonal group $O(V)$, the group of (unipotent) upper triangular $n\times n$ matrices acting naturally on $\field^n$, the group of monomial $n\times n$ matrices acting on $\field^n$, and direct products of such linear groups all satisfy the conditions of 
Proposition~\ref{prop:pseudo-refl-indep}. }
\end{example}

\begin{center} {\bf Acknowledgements}\end{center} 

I thank Endre Szab\'o and Tam\'as Szamuely for helpful discussions, as well as an anonymous referee of a first version of this 
note for some references, that helped me to gain far 
more conclusive final results from the basic ideas.


 \end{document}